\documentclass[11pt]{amsart} 
\usepackage{amssymb,amsmath,latexsym,enumerate,graphicx,bbm,mathptmx,microtype}
\newtheorem{theorem}{Theorem} 

\newtheorem{lemma}[theorem]{Lemma}

\newtheorem*{rem}{Remarks}

\newtheorem*{exam}{Examples}

\newcommand\commentout[1]{}
\renewcommand\emptyset{\varnothing}
\def\conv{\operatorname{conv}}
\newcommand\ehr{\operatorname{ehr}}
\newcommand\Ehr{\operatorname{Ehr}}
\newcommand\cone{\operatorname{cone}}
\newcommand\mult{\operatorname{mult}}
\newcommand\Des{\operatorname{Des}}
\newcommand\Asc{\operatorname{Asc}}

\newcommand\maj{\operatorname{maj}}
\newcommand\amaj{\operatorname{amaj}}
\def\th{^{\text{th}}}
\def\Kappa{\mathrm{K}}
\def\Z{\mathbbm{Z}}
\def\Q{{\mathcal Q}}
\def\R{\mathbbm{R}}
\def\C{\mathbbm{C}}
\def\V{\mathcal{V}}

\def\P{\mathcal{P}}
\def\F{\mathcal{F}}
\def\G{\mathcal{G}}
\def\H{\mathcal{H}}

\def\m{{\mathbf m}}
\def\v{{\mathbf v}}
\def\x{{\mathbf x}}
\def\0{{\mathbf 0}}
\def\1{{\mathbf 1}}

\def\chrom{c}
\def\h{h}

\title{Combinatorial Reciprocity Theorems}

\author{Matthias Beck}
\address{Department of Mathematics\\
         San Francisco State University\\
         San Francisco, CA 94132\\
         U.S.A.}
\email{mattbeck@sfsu.edu}

\subjclass[2000]{05A15, 05C15, 05C31, 11H06, 52C07, 52C35}
\keywords{Combinatorial reciprocity theorem, rational generating function, convex polyhedron, Euler--Poincar\'e relation, hyperplane arrangement, lattice point, lattice polytope, Ehrhart polynomial, chromatic polynomial, acyclic orientation of a graph, inside-out polytope, poset, $P$-partition, permutation statistics}
\thanks{\emph{Danke} to J\"org Rambau for initiating this survey article and to an anonymous referee for helpful suggestions.
The author was partially supported by the NSF (DMS-0810105).
This survey paper is a (vastly) compressed version of an upcoming book with the same title; a preprint can be found at {\tt math.sfsu.edu/beck/crt.html}.}
\date{30 December 2011. To appear in \emph{Jahresbericht der DMV}}

\begin{document}

\begin{abstract}
A common theme of enumerative combinatorics is formed by counting functions that are polynomials evaluated at positive integers. 
In this expository paper, we focus on four families of such counting functions connected to hyperplane arrangements, lattice points in polyhedra, proper colorings of graphs, and $P$-partitions.
We will see that in each instance we get interesting information out of a counting function when we evaluate it at a \emph{negative} integer (and so, a priori the counting function does not make sense at this number). Our goals are to convey some of the charm these ``alternative" evaluations of counting functions exhibit, and to weave a unifying thread through various combinatorial reciprocity theorems by looking at them through the lens of geometry, which will include some scenic detours through other combinatorial concepts.
\end{abstract}

\dedicatory{Dedicated to my friend and mentor Tom Zaslavsky}

\maketitle



\section{Introduction}

A common theme of enumerative combinatorics is formed by counting functions that are polynomials evaluated at positive integers. 
To be as concrete as possible, we focus on four families of such counting functions. 
We will see that in each instant we get interesting information out of a counting function when we evaluate it at a \emph{negative} integer (and so, a priori the counting function does not make sense at this number). Our goals are to convey some of the charm these ``alternative" evaluations of counting functions exhibit, and to weave a unifying thread through various combinatorial reciprocity theorems by looking at them through the lens of geometry, which will include some scenic detours through other combinatorial concepts.

We have tried to keep this expository paper self contained, requiring only a few basic, well-known facts about polyhedra (such as the Euler--Poincar\'e relation) and a healthy dose of enthusiasm for exercises, which we have implicitly spread throughout these notes.
We start by introducing the main players of this story.


\subsection{Hyperplane Arrangements}
A \emph{hyperplane arrangement} $\H$ is a finite collection of hyperplanes in $\R^d$.
A \emph{flat} of $\H$ is a nonempty intersection of some of the hyperplanes in $\H$; we always include $\R^d$ among the flats\footnote{$\R^d$ is the flat you obtain when you don't intersect anything.}.
Flats are naturally ordered by (reverse) set inclusion; see Figure \ref{hyparrposet3fig} for an example.
A \emph{region} of $\H$ is a maximal connected component of $\R^d \setminus \bigcup \H$.
Our first goal is to count the regions of a hyperplane arrangement~$\H$. 

\def\JPicScale{0.6}
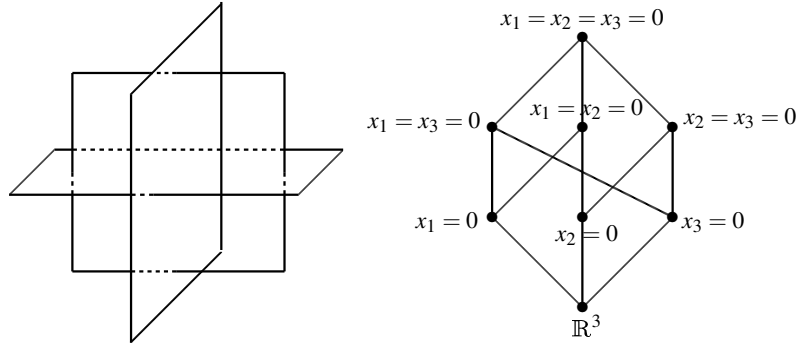
\begin{figure}[htb]
\begin{center}
\ifx\JPicScale\undefined\def\JPicScale{1}\fi
\unitlength \JPicScale mm
\begin{picture}(165,87.75)(0,0)
\linethickness{0.3mm}
\put(40,72){\line(1,0){24}}
\linethickness{0.3mm}
\put(64,50){\line(0,1){22}}
\linethickness{0.3mm}
\put(17,50){\line(0,1){22}}
\linethickness{0.3mm}
\put(17,72){\line(1,0){17}}
\linethickness{0.3mm}
\multiput(34,72)(1.71,0){4}{\line(1,0){0.86}}
\linethickness{0.3mm}
\put(40,28){\line(1,0){24}}
\linethickness{0.3mm}
\put(17,28){\line(1,0){13}}
\linethickness{0.3mm}
\multiput(30,28)(1.82,0){6}{\line(1,0){0.91}}
\linethickness{0.3mm}
\multiput(3,45)(0.12,0.12){83}{\line(1,0){0.12}}
\linethickness{0.3mm}

\linethickness{0.3mm}
\put(30,12.25){\line(0,1){55}}
\multiput(30,12.25)(0.12,0.12){167}{\line(0,1){0.12}}
\multiput(30,67.25)(0.12,0.12){167}{\line(0,1){0.12}}
\put(50,32.75){\line(0,1){55}}
\linethickness{0.3mm}
\multiput(67,45)(0.12,0.12){83}{\line(1,0){0.12}}
\linethickness{0.3mm}
\put(13,55){\line(1,0){4}}
\linethickness{0.3mm}
\multiput(17,55)(2,0){24}{\line(1,0){1}}
\linethickness{0.3mm}
\put(64,55){\line(1,0){13}}
\linethickness{0.3mm}
\put(3,45){\line(1,0){27}}
\linethickness{0.3mm}
\put(35,45){\line(1,0){32}}
\linethickness{0.3mm}
\multiput(30,45)(2,0){3}{\line(1,0){1}}
\linethickness{0.3mm}
\multiput(17,45)(0,2){3}{\line(0,1){1}}
\linethickness{0.3mm}
\put(17,28){\line(0,1){17}}
\linethickness{0.3mm}
\multiput(64,45)(0,2){3}{\line(0,1){1}}
\linethickness{0.3mm}
\put(64,28){\line(0,1){17}}
\linethickness{0.3mm}
\put(130,80){\circle*{2.5}}

\linethickness{0.3mm}
\put(130,60){\circle*{2.5}}

\linethickness{0.3mm}
\put(110,60){\circle*{2.5}}

\linethickness{0.3mm}
\put(110,40){\circle*{2.5}}

\linethickness{0.3mm}
\put(130,40){\circle*{2.5}}

\linethickness{0.3mm}
\put(130,20){\circle*{2.5}}

\linethickness{0.3mm}
\put(150,60){\circle*{2.5}}

\linethickness{0.3mm}
\put(150,40){\circle*{2.5}}

\linethickness{0.3mm}
\multiput(110,40)(0.12,-0.12){167}{\line(1,0){0.12}}
\linethickness{0.3mm}
\put(130,20){\line(0,1){60}}
\linethickness{0.3mm}
\put(110,40){\line(0,1){20}}
\linethickness{0.3mm}
\multiput(110,40)(0.12,0.12){167}{\line(1,0){0.12}}
\linethickness{0.3mm}
\multiput(130,40)(0.12,0.12){167}{\line(1,0){0.12}}
\linethickness{0.3mm}
\multiput(110,60)(0.24,-0.12){167}{\line(1,0){0.24}}
\linethickness{0.3mm}
\multiput(110,60)(0.12,0.12){167}{\line(1,0){0.12}}
\linethickness{0.3mm}
\multiput(130,80)(0.12,-0.12){167}{\line(1,0){0.12}}
\linethickness{0.3mm}
\put(150,40){\line(0,1){20}}
\linethickness{0.3mm}
\multiput(130,20)(0.12,0.12){167}{\line(1,0){0.12}}
\put(132,83){\makebox(0,0)[cc]{}}

\put(132,83){\makebox(0,0)[cc]{}}

\put(129,84){\makebox(0,0)[cc]{}}

\put(129,84){\makebox(0,0)[cc]{}}

\put(98,78){\makebox(0,0)[cc]{}}

\put(100,39){\makebox(0,0)[cc]{\footnotesize $x_1=0$}}

\put(95,61){\makebox(0,0)[cc]{\footnotesize $x_1=x_3=0$}}

\put(95,40){\makebox(0,0)[cc]{}}

\put(131,16){\makebox(0,0)[cc]{\footnotesize $\R^3$}}

\put(131,36){\makebox(0,0)[cc]{\footnotesize $x_2=0$}}

\put(159,39){\makebox(0,0)[cc]{\footnotesize $x_3=0$}}

\put(165,62){\makebox(0,0)[cc]{\footnotesize $x_2=x_3=0$}}

\put(131,64){\makebox(0,0)[cc]{\footnotesize $x_1=x_2=0$}}

\put(130,84){\makebox(0,0)[cc]{\footnotesize $x_1=x_2=x_3=0$}}

\end{picture}
\end{center}
\caption{An arrangement of three coordinate hyperplanes and its flats.}\label{hyparrposet3fig}
\end{figure}

The \emph{M\"obius function} of $\H$ is defined on the set of all flats of $\H$ recursively through
\begin{equation}\label{mobiusdef}
  \mu (F) =
  \begin{cases}
    1 & \text{ if } F = \R^d , \\
    \displaystyle - \sum_{ G \supsetneq F } \mu (G) & \text{ otherwise. }
  \end{cases}
\end{equation}
(M\"obius functions can be defined in much greater generality, as we will see in Section \ref{hyperplsection}.)
The M\"obius function, in turn, allows us to define the \emph{characteristic polynomial} of $\H$ by
\[
  \h_\H (t) := \sum_{ F \in L(\H) } \mu(F) \, t^{ \dim F } .
\]
Here are some examples of classic families of hyperplane arrangements and their characteristic polynomials, whose computation makes for a fun exercise.
\begin{itemize}
\item For the \emph{Boolean arrangement} $\H = \left\{ x_j = 0 : \, 1 \le j \le d \right\}$, $\h_\H (t) = (t-1)^d$.
\item For the \emph{braid arrangement} $\H = \left\{ x_j = x_k : \, 1 \le j < k \le d \right\}$, $\h_\H (t) = t (t-1) (t-2) \cdots (t-d+1)$.
\item For an arrangement $\H$ in $\R^d$ consisting of $n$ hyperplanes in general position,
\[
  \h_\H (t) = \binom n 0 t^d - \binom n 1 \, t^{ d-1 } + \binom n 2 t^{ d-2 } - \dots + (-1)^d \binom n d \, .
\]
\end{itemize}
The astute reader will notice that each of these characteristic polynomials bear the number of regions of the hyperplane arrangement in question as the special evaluation $\left| \h_\H(-1) \right|$. For example, the braid arrangement dissects $\R^d$ into $d!$ regions.
This is not an accident:

\begin{theorem}[Zaslavsky \cite{zaslavskythesis}]\label{zaslavskythm}
Suppose $\H$ is a hyperplane arrangement in $\R^d$. Then $(-1)^{ d } \, \h_\H (-1)$ equals the number of regions of $\H$.
\end{theorem}


\subsection{Ehrhart Polynomials}

A \emph{lattice polytope}
is the convex hull (in $\R^d$) of finitely many points in $\Z^d$.
For such a polytope $\P$, we define
\[
  \ehr_\P (t) := \# \left( t \P \cap \Z^d \right) ,
\]
the number of integer lattice points in the $t\th$ dilate of $\P$, where $t$ is a positive integer.
As an example, consider the triangle $\Delta \subset \R^2$ with vertices $(0,0)$, $(1,0)$, and $(0,1)$. It comes with the lattice-point enumerator
\[
  \ehr_\Delta (t) = \# \left\{ (m,n) \in \Z^2 : \, m, n \ge 0 , \, m+n \le t \right\} .
\]

\begin{figure}[htb]
\begin{center}
\includegraphics[totalheight=1.5in]{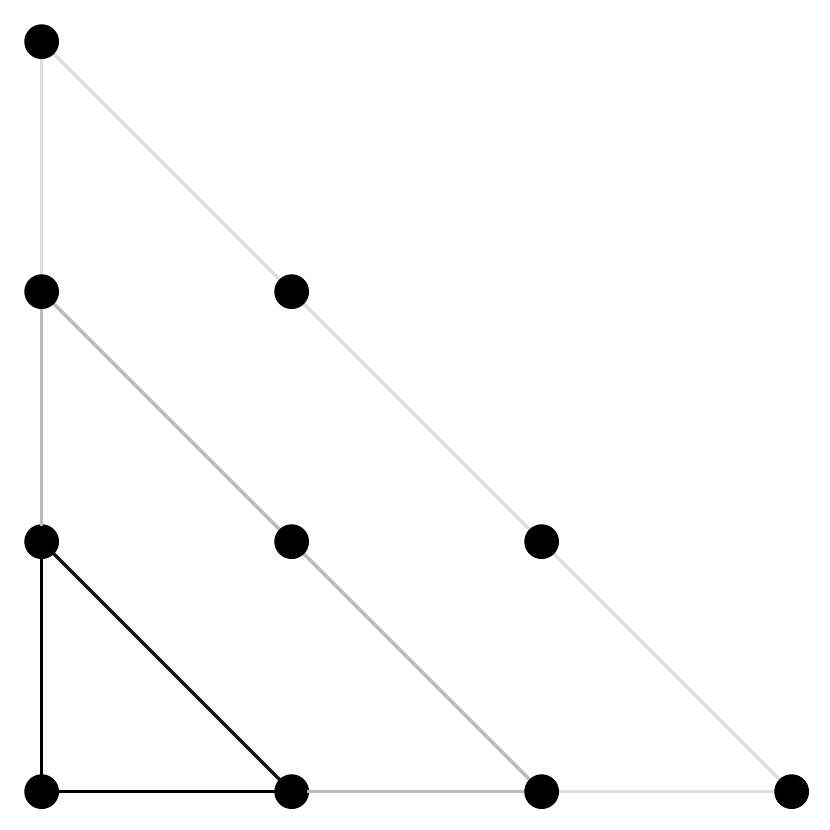}
\end{center}
\caption{A sample lattice-point problem.}\label{introtrianglefig}
\end{figure}

\noindent
A moment's thought (or a look at Figure \ref{introtrianglefig}) reveals that $\ehr_\Delta (t)$ is given by \emph{triangular numbers}:
\[
  \ehr_\Delta (t) = \binom{ t+2 }{ 2 } = \frac 1 2 (t+1)(t+2) \, .
\]
If we evaluate this polynomial at $-t$, we obtain
\[
  \ehr_\Delta (-t) = \frac 1 2 (t-1)(t-2) = \binom{ t-1 }{ 2 } \, ,
\]
which happens to be the function enumerating the \emph{interior} lattice point in $t \Delta$, by another counting argument for triangular numbers (just draw a picture of the interior lattice points in $t \Delta$!).
So in our example, we obtain the functional relation
\[
  \ehr_\Delta (-t) = \ehr_{\Delta^\circ} (t) \, ,
\]
where $\Delta^\circ$ denotes the interior of $\Delta$.
For example, the evaluations $\ehr_\Delta (-1) = \ehr_\Delta (-2) = 0$ point to the fact that neither $\Delta$ nor $2 \Delta$ contain any interior lattice points.
Once more this is far from accidental:

\begin{theorem}[Ehrhart \cite{ehrhartpolynomial}, Macdonald \cite{macdonald}]\label{ehrmacthm}
If $\P$ is a lattice $d$-polytope, then for positive integers $t$, the counting function $\ehr_\P(t)$ is a polynomial in $t$.
When this polynomial is evaluated at negative integers, we obtain
\[
  \ehr_\P (-t) = (-1)^d \, \ehr_{ \P^\circ } (t) \, .
\]
\end{theorem}


\subsection{Chromatic Polynomials}
\begin{figure}[htb]
\begin{center}
\includegraphics[totalheight=.8in]{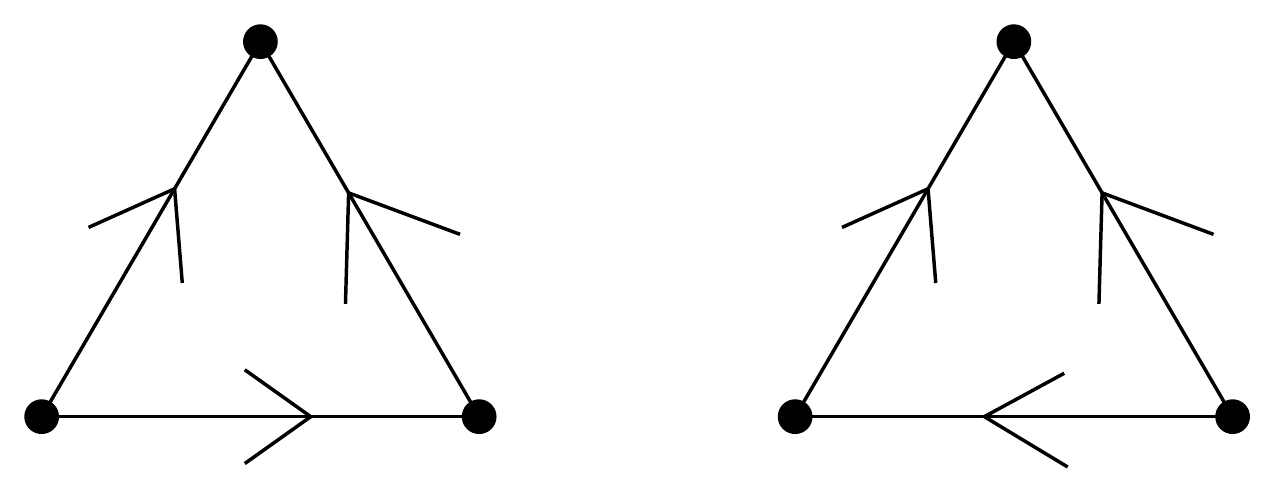}
\end{center}
\caption{Two of the acyclic orientations of $K_3$.}\label{introacyclicfig}
\end{figure}

Let $G = (V, E)$ be a graph.
The \emph{chromatic polynomial} $\chrom_G(t)$ (whose roots can be traced to Birkhoff \cite{birkhoffcoloring} and Whitney \cite{whitneychromatic}) is the counting function that enumerates all \emph{proper $t$-colorings}, i.e., labellings $\x \in [t]^V$ such that adjacent nodes get different labels: $ij \in E \ \Longrightarrow \ x_i \ne x_j$. (Here $[t] := \{ 1, 2, \dots, t \}$.)
For example, the graph $K_3$ with three nodes, any pair of which is adjacent, has chromatic polynomial
\[
  \chrom_{ K_3 } (t) = t (t-1) (t-2) \, ,
\]
as all three nodes get different labels.
When we evaluate this chromatic polynomial at $-1$, we obtain
\[
  \chrom_{ K_3 } (-1) = -6 \, ,
\]
which is, up to a sign, the number of \emph{acyclic orientations} of $K_3$, namely, those orientations that do not contain any coherently oriented cycle (see Figure \ref{introacyclicfig}). 
The evaluation $\chrom_{ K_3 } (-1)$ is part of a much more general phenomenon, for which we need one more definition:
An orientation of $G$ and a (not necessarily proper) $t$-coloring $\x \in [t]^V$ are \emph{compatible} if $x_j \ge x_i$ whenever there is an edge oriented from $i$ to $j$.

\begin{theorem}[Stanley \cite{stanleyacyclic}]\label{coloringreciprocitythm}
Let $G = (V, E)$ be a graph with finite node set $V$. Then $(-1)^{|V|} \, \chrom_G(-t)$ equals the number of pairs consisting of a $t$-coloring and a compatible acyclic orientation of $G$.
In particular, $(-1)^{|V|} \, \chrom_G(-1)$ counts all acyclic orientations of~$G$.
\end{theorem}


\subsection{P-partitions}\label{ppartintrosection}

Our final example originates in the world of integer partitions, with a connection to partially-ordered sets (\emph{posets}).
Recall that a \emph{partition} of the integer $t$ is a sequence $\left( x_1, x_2, \dots, x_d \right)$ of nonnegative integers such that
\begin{equation}\label{ppartitiondef}
  t = x_1 + x_2 + \dots + x_d
  \qquad \text{ and } \qquad
  x_1 \ge x_2 \ge \dots \ge x_d \, .
\end{equation}
There are instances when we are interested in writing an integer $t$ in the form
\begin{equation}\label{compositiondef}
  t = x_1 + x_2 + \dots + x_d \, ,
\end{equation}
i.e., \emph{without} the restriction $x_1 \ge x_2 \ge \dots \ge x_d$; then we call $\left( x_1, x_2, \dots, x_d \right) \in \Z_{ \ge 0 }^d$ a \emph{composition} of $t$.
The theory of $P$-partitions\footnote{here $P$ stands for a specific poset---for which we tend to use greek letters such as $\Pi$ to avoid confusions with polytopes.} allows us to interpolate between \eqref{ppartitiondef} and \eqref{compositiondef}; that is, we will study compositions of $t$ that satisfy \emph{some} of the inequalities (implied by) $x_1 \ge x_2 \ge \dots \ge x_d$. A natural way to introduce such a subset of inequalities is through a poset $\Pi$, whose relation we denote by $\preceq$.
A \emph{$\Pi$-partition} is an order-reversing map $\x: \Pi \to \Z_{ \ge 0 }$, i.e.,
\[
  a \preceq b
  \qquad \Longrightarrow \qquad
  x_a \ge x_b \, .
\]
A \emph{strict $\Pi$-partition} is a map $\x: \Pi \to \Z_{ \ge 0 }$ such that
\[
  a \prec b
  \qquad \Longrightarrow \qquad
  x_a > x_b \, .
\]
In either case, if $\sum_{ a \in \Pi } x_a = t$ then we call $\x$ a \emph{(strict) $\Pi$-partition of $t$}.
Let $p_\Pi(t)$ denote the number of $\Pi$-partitions of $t$, with generating function
\[
  P_\Pi(z)
  := \sum_{ t \ge 0 } p_\Pi(t) \, z^t
  = \sum_{ \x } z^{ x_1 + x_2 + \dots + x_d } ,
\]
where the last sum is taken over all $\Pi$-partitions.
Analogously, we define the number of \emph{strict} $\Pi$-partitions of $t$ as $p_\Pi^\circ(t)$, with accompanying generating function $P_\Pi^\circ(z)$.

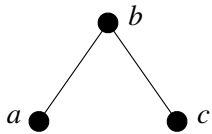
\begin{figure}[htb]
\begin{center}
\begin{picture}(70,50)
\includegraphics[totalheight=.7in]{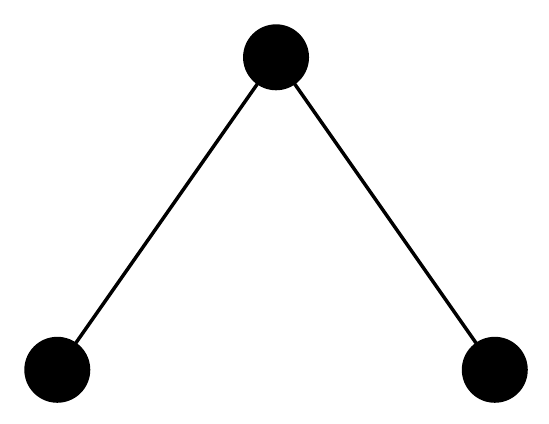}
\end{picture}
\begin{picture}(0,0)
  \put(-78,6){$a$}
  \put(-32,43){$b$}
  \put(-6,6){$c$}
\end{picture}
\end{center}
\caption{A sample poset.}\label{introposetfig}
\end{figure}

Here are three basic examples, which we invite the reader to work out, and which illustrate the way $P$-partitions are situated between partitions and compositions:
\begin{enumerate}[(i)]
\item If $\Pi$ is a \emph{chain} with $d$ elements (whose elements are totally ordered), a $\Pi$-partition of $t$ is a partition of $t$ in the sense of~\eqref{ppartitiondef}, with generating functions
\[
  P_\Pi(z) = \frac{ 1 }{ (1-z) (1-z^2) \cdots (1-z^d) }
  \qquad \text{ and } \qquad
  P_\Pi^\circ(z) = \frac{ z^{ \binom d 2 } }{ (1-z) (1-z^2) \cdots (1-z^d) } \, .
\]
\item If $\Pi$ is an \emph{antichain} with $d$ elements (whose elements have no relation whatsoever), a $\Pi$-partition of $t$ is a composition of $t$ in the sense of~\eqref{compositiondef}, with generating functions
\[
  P_\Pi(z) = \frac{ 1 }{ (1-z)^d } = P_\Pi^\circ(z) \, .
\]
\item If $\Pi$ be the poset pictured in Figure~\ref{introposetfig}, then
\[
  P_\Pi(z) = \frac{ 1 }{ (1-z)^2 (1-z^3) }
  \qquad \text{ and } \qquad
  P_\Pi(z) = \frac{ z^2 }{ (1-z)^2 (1-z^3) } \, .
\]
\end{enumerate}

By now it should come as no surprise that there is a combinatorial reciprocity theorem relating these generating functions. Since reciprocity for a (quasi-)polynomial\footnote{In general, $p_\Pi(t)$ is a \emph{quasipolynomial}, a term we will define in Section \ref{ehrhartgeneralsection}.} means replacing the variable $t$ by $-t$, when we express reciprocity in terms of generating functions, we should replace the variable $z$ by~$\frac 1 z$.

\begin{theorem}[Stanley \cite{stanleythesis}]\label{pipartitionreciprocitythm}
Given a finite poset $\Pi$, the rational functions $P_\Pi(z)$ and $P_\Pi^\circ(z)$ are related by
\[
  P_\Pi ( \tfrac 1 z ) = (-z)^{ |\Pi| } P_\Pi^\circ(z) \, .
\]
\end{theorem}


Let us reiterate the common thread that can be weaved through Theorems \ref{zaslavskythm}--\ref{pipartitionreciprocitythm}. 
Each of them is an instance of a \emph{combinatorial reciprocity theorem}: a combinatorial function, which is a priori defined on the positive integers,
\begin{enumerate}[(i)]
\item can be algebraically extended beyond the positive integers (e.g., because it is a polynomial), and
\item has (possibly quite different) meaning when evaluated at negative integers.
\end{enumerate}

We will illustrate a geometric approach to the above reciprocity theorems, by mixing lattice points, polyhedra, and hyperplane arrangements, in the sense that we interpret the objects we'd like to count as lattice points, subject to some linear constraints (giving rise to a polyhedron), with an interplay given by further linear conditions (giving rise to hyperplanes).

Thus Theorems \ref{zaslavskythm} and \ref{ehrmacthm} can be used as building blocks to prove Theorems \ref{coloringreciprocitythm} and \ref{pipartitionreciprocitythm} (though this is not historically how the original proofs surfaced).
We will give the main ideas for proofs of Theorems \ref{zaslavskythm} and \ref{ehrmacthm} in Sections \ref{hyperplsection} and \ref{ehrhartsection}, respectively.
Both proofs rely on variants of the Euler--Poincar\'e relation of a polyhedron, which in itself can be thought of as a combinatorial reciprocity theorem.
In Sections \ref{coloringiopsection} and \ref{iopsection}, we give a proof of Theorem \ref{coloringreciprocitythm} to illustrate how the introduction of ``forbidden" hyperplanes into Ehrhart's theory of lattice-point enumeration in polytopes allows us to prove old and new combinatorial reciprocity theorems geometrically.
A second way of arranging hyperplanes with polyhedra, as triangulation hyperplanes, is illustrated in Section \ref{ppartsection}, which contains a proof of Theorem \ref{pipartitionreciprocitythm} and connections to permutation statistics.


\section{The Euler--Poincar\'e Relation and Zaslavsky's Theorem}\label{hyperplsection}

\subsection{Polyhedra}
A \emph{(convex) polyhedron} $\P$ is the intersection of finitely many (affine) halfspaces in $\R^d$.
Bounded polyhedra are \emph{polytopes}; the fact that they can also be described as the convex hull of finitely many points in $\R^d$ is the famous (and nontrivial) \emph{Minkowski--Weyl Theorem} (see., e.g., \cite[Lecture 1]{ziegler}).
An even more famous theorem concerns the polynomial
\[
  f_\P (t) = \sum_\F t^{ \dim \F } ,
\]
where we sum over all (nonempty) faces\footnote{A \emph{face} of $\P$ is a set of the form $\P \cap H$, where $H$ is a hyperplane that bounds a half space containing $\P$; we always include $\P$ itself (and sometimes $\emptyset$) in the list of faces of $\P$.} of $\P$:

\begin{theorem}[Euler \cite{euler2,euler1}, Poincar\'e \cite{poincare}] \label{eulerrelation}
Suppose $\P = \V + \Q$ is a polyhedron, where $\V$ is a vector space and $\Q$ is a polyhedron that contains no lines.\footnote{Here $+$ refers to Minkowski (point-wise) sum; it is an easy fact that every polyhedron can be written as a sum of a vector space and a polyhedron that contains no lines.} Then
\[
  f_\P (-1) = 
  \begin{cases}
    (-1)^{ \dim \V } & \text{ if $\Q$ is bounded,} \\
    0 & \text{ if $\Q$ is unbounded.}
  \end{cases}
\]
\end{theorem}

Our formulation of the Euler--Poincar\'e relation suggests that it can be viewed as a combinatorial reciprocity theorem in its own right (and in a sense all other such reciprocity theorems are based on it).
The number $f_\P (-1)$ is usually called the \emph{Euler characteristic} of~$\P$.

The function $\mu(F)$ we defined in \eqref{mobiusdef} is a special case of the following construct.
For a general poset $\Pi$ equipped with a relation $\preceq$, we define its \emph{M\"obius function} recursively through
\begin{equation}\label{genmobiusdef}
  \mu (x,y) =
  \begin{cases}
    0 & \text{ if } x \not\preceq y \, , \\
    1 & \text{ if } x=y \, , \\
    \displaystyle - \sum_{ x \preceq z \prec y } \mu (x,z) & \text{ if } x \prec y \, .
  \end{cases}
\end{equation}
The central result for these functions, which is a fun exercise, is \emph{M\"obius inversion}: for $f, g \in \C^\Pi$,
\begin{equation}\label{mobiusinversion}
  f(x) = \sum_{ y \succeq x } g(y) \qquad \Longleftrightarrow \qquad g(x) = \sum_{ y \succeq x } \mu(x,y) \, f(y) \, .
\end{equation}
The M\"obius function of a poset gives rise to a generalization of the \emph{inclusion--exclusion principle} (which follows from M\"obius inversion for the poset of intersections of a given family of sets); see, e.g., \cite[Chapter 3]{stanleyec1} for much more about M\"obius functions.

One can view the Euler--Poincar\'e relation (Theorem \ref{eulerrelation}) in the light of the M\"obius function of the poset formed by all faces of a polytope $\P$: By the recursive definition \eqref{genmobiusdef} of $\mu$, we have for any nonempty face $\F$,
\begin{equation}\label{eulerrevisited1}
  \sum_{ \emptyset \subseteq \G \subseteq \F } \mu(\emptyset, \G) = 0 \, ,
\end{equation}
where we sum over all faces $\G$ of $\F$, including $\emptyset$ and $\F$ itself.
On the other hand, each face $\F$ is again a polytope, and so the Euler--Poincar\'e relation (Theorem \ref{eulerrelation}) says that
\[
  \sum_{ \emptyset \subsetneq \G \subseteq \F } (-1)^{ \dim \G } = 1 \, .
\]
But this implies (if we add the empty face to our sum, giving it dimension $-1$) that
\[
  \sum_{ \emptyset \subseteq \G \subseteq \F } (-1)^{ \dim \G } = 0 \, .
\]
Since both this equation and \eqref{eulerrevisited1} hold for \emph{any} face $\F$ and
$
  \mu(\emptyset, \emptyset) = 1 = (-1)^{ \dim \emptyset + 1 } ,
$
we recursively compute
\[
  \mu(\emptyset, \F) = (-1)^{ \dim \F + 1 }
\]
for all faces $\F \subseteq \P$; more generally, one can show
$
  \mu(\G,\F) = (-1)^{ \dim \F - \dim \G } .
$


\subsection{Hyperplane arrangements}
Now we connect the above concepts to a hyperplane arrangement $\H$ in $\R^d$.
The flats of $\H$ form a poset $L(\H)$ which we order by reverse set inclusion:
\[
  F \preceq G \quad \Longleftrightarrow \quad F \supseteq G \, .
\]
Thus the M\"obius function $\mu(F)$ we defined in \eqref{mobiusdef} equals the special evaluation $\mu \left( \R^d, F \right)$ of the M\"obius function of $L(\H)$.

A face of any of the regions of $\H$ is called a \emph{face} of $\H$.
Given a flat $F$ of $\H$, we can create the hyperplane arrangement \emph{induced} by $\H$ on $F$, namely,
\[
  \H^F := \left\{ H \cap F : \, H \in \H, \, H \cap F \ne \emptyset \right\} .
\]
The proof of Zaslavsky's Theorem \ref{zaslavskythm} is based on the observation that
each face $f$ of $\H^F$ is a region of $\H^G$ for some flat $G \subseteq F$ (more precisely, $G$ is the affine span of $f$), and so
\[
  \sum_{ f \text{ face of } \H^F } (-1)^{ \dim f } = \sum_{ G \subseteq F } (-1)^{ \dim G } \, r(\H^G) \, .
\]
But the left-hand side is simply the Euler characteristic of $F$, which is $(-1)^{ \dim F }$ (by Theorem \ref{eulerrelation}).
Thus
\[
  (-1)^{ \dim F } = \sum_{ G \subseteq F } (-1)^{ \dim G } \, r(\H^G) \, ,
\]
and we can use M\"obius inversion \eqref{mobiusinversion}:
\[
  (-1)^{ \dim F } r(\H^F) = \sum_{ G \subseteq F } \mu(F, G) \, (-1)^{ \dim G } .
\]
For $F = \R^d$ this gives Theorem~\ref{zaslavskythm}:
\[
  (-1)^d \, r(\H) = \sum_{ G \in L(\H) } \mu(G) \, (-1)^{ \dim G } = \h_\H (-1) \, .
\]

For more on the combinatorics of hyperplane arrangements, we recommend the survey article \cite{stanleyparkcityhyperplanes}.
For numerous interesting topological considerations that arise from the study of hyperplane arrangements over $\C$, see~\cite{orlikterao}.


\section{Ehrhart--Macdonald Reciprocity}\label{ehrhartsection}

\subsection{Lattice Simplices}\label{ehrhartsimplexsection}

In this section, we will give an idea why Theorem \ref{ehrmacthm} is true.
We will first show how to prove it for lattice \emph{simplices}, each of which is the convex hull of $d+1$ affinely independent points in $\R^n$ (and in this section we will assume $n=d$).
We form the \emph{cone over} such a simplex $\Delta$
\[
  \cone(\Delta) := \sum_{ \v \text{ vertex of } \Delta } \R_{ \ge 0 } (\v, 1)
\]
by lifting the vertices of $\Delta$ into $\R^{ d+1 }$ onto the hyperplane $x_{ d+1 } = 1$ and taking the nonnegative span of this ``lifted version" of $\Delta$; see Figure \ref{coneoversimplexfig} for an illustration.
\def\JPicScale{0.7}
\begin{figure}[htb]
\begin{center}
\ifx\JPicScale\undefined\def\JPicScale{1}\fi
\unitlength \JPicScale mm
\begin{picture}(160,30)(0,0)
\linethickness{0.3mm}
\put(45,10){\circle*{2.5}}

\linethickness{0.3mm}
\put(0,0){\line(1,0){70}}
\put(70,0){\vector(1,0){0.12}}
\linethickness{0.9mm}
\put(15,10){\line(1,0){30}}
\linethickness{0.3mm}
\put(25,0){\line(0,1){30}}
\put(25,30){\vector(0,1){0.12}}
\linethickness{0.3mm}
\put(15,10){\circle*{2.5}}

\linethickness{0.3mm}
\multiput(25,0)(0.24,0.12){188}{\line(1,0){0.24}}
\linethickness{0.3mm}
\multiput(0,25)(0.12,-0.12){208}{\line(1,0){0.12}}
\linethickness{0.3mm}
\put(135,10){\circle*{2.5}}

\linethickness{0.3mm}
\put(90,0){\line(1,0){70}}
\put(160,0){\vector(1,0){0.12}}
\linethickness{0.3mm}
\put(115,0){\line(0,1){30}}
\put(115,30){\vector(0,1){0.12}}
\linethickness{0.3mm}
\put(105,10){\circle*{2.5}}

\linethickness{0.3mm}
\multiput(115,0)(0.24,0.12){188}{\line(1,0){0.24}}
\linethickness{0.3mm}
\multiput(90,25)(0.12,-0.12){208}{\line(1,0){0.12}}
\linethickness{0.3mm}
\multiput(105,10)(1.74,0.87){12}{\multiput(0,0)(0.22,0.11){4}{\line(1,0){0.22}}}
\linethickness{0.3mm}
\multiput(125,20)(1.33,-1.33){8}{\multiput(0,0)(0.11,-0.11){6}{\line(1,0){0.11}}}
\put(54,9){\makebox(0,0)[cc]{$(2,1)$}}

\put(55,7.5){\makebox(0,0)[cc]{}}

\put(5,9){\makebox(0,0)[cc]{$(-1,1)$}}

\end{picture}
\end{center}
\caption{The cone over the one-dimensional simplex $[-1,2]$ and its fundamental parallelepiped.}\label{coneoversimplexfig}
\end{figure}
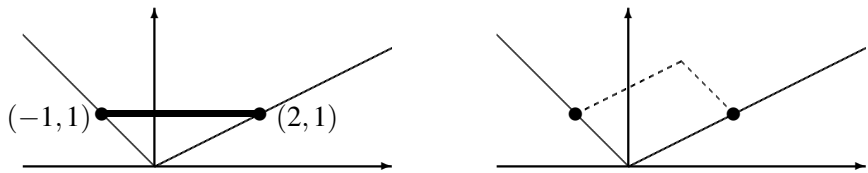
The reason for coning over $\Delta$ is that we can see a copy of the dilate $t \Delta$ as the intersection of $\cone(\Delta)$ with the hyperplane $x_{ d+1 } = t$; we will say that theses points are \emph{at height} $t$. So the \emph{Ehrhart series}
\begin{equation}\label{ehrhartseriesdef}
  \Ehr_\Delta(z) := 1 + \sum_{ t > 0 } \ehr_\Delta(t) \, z^t
\end{equation}
can be computed through
\[
  \Ehr_\Delta(z) = \sum_{ t \ge 0 } \# \left( \text{lattice points in } \cone(\Delta) \text{ at height } t \right) z^t .
\]
We use a tiling argument to compute this generating function. Namely, let
\[
  \Q := \sum_{ \v \text{ vertex of } \Delta } [0,1) (\v, 1) \, ,
\]
the \emph{fundamental parallelepiped} of $\cone(\Delta)$. Then we can tile $\cone(\Delta)$ by translates of~$\Q$:
\[
  \cone(\Delta) = \bigcup_{ \m \in \Z_{ \ge 0 }^{ d+1 } } \left( \sum_{ \v \text{ vertex of } \Delta } \!\!\!\!m_\v (\v,1) \ +\  \Q \right) ,
\]
and this union is disjoint (because $\Q$ is half open).
Every lattice point in $\cone(\Delta)$ is a translate of such a nonnegative integral combination of the $(\v,1)$'s by a lattice point in $\Q$ (and this representation is unique). Translated into generating-function language, this gives
\[
  \Ehr_\Delta(z)
  = \left( \frac{ 1 }{ 1-z } \right)^{ d+1 } \sum_{ t \ge 0 } \# \left( \text{lattice points in } \Q \text{ at height } t \right) z^t .
\]
The sum on the right is a polynomial $h(z)$ of degree at most $d$, and it is a basic exercise to deduce from the rational-function form of $\Ehr_\Delta(z)$ that $\ehr_\Delta(t)$ is a polynomial.
This proves the first part of Theorem \ref{ehrmacthm} in the simplex case.
Towards the second part, we compute
\[
  \Ehr_\Delta \left( \frac 1 z \right)
  = \sum_{ t \ge 0 } \ehr_\Delta(t) \, z^{ -t }
  = \frac{ h \left( \frac 1 z \right) }{ \left( 1 - \frac 1 z \right)^{ d+1 } } 
  = (-1)^{ d+1 } \frac{ z^{ d+1 } h \left( \frac 1 z \right) }{ (1-z)^{ d+1 } } 
\]
and so by an easy exercise about generating functions,
\[
  \sum_{ t < 0 } \ehr_\Delta(t) \, z^{ -t }
  = \sum_{ t>0 } \ehr_\Delta(-t) \, z^t
  = (-1)^d \frac{ z^{ d+1 } h \left( \frac 1 z \right) }{ (1-z)^{ d+1 } } \, .
\]
Inspired by this, we define \label{ehrhartseriesinteriornotation}
\begin{equation}\label{ehrhartseriesinteriordef}
  \Ehr_{ \Delta^\circ } (z) := \sum_{ t > 0 } \ehr_{ \Delta^\circ } (t) \, z^t ,
\end{equation}
and so proving the reciprocity theorem $\ehr_\Delta (-t) = (-1)^d \, \ehr_{ \Delta^\circ } (t)$ is equivalent to proving
\begin{equation}\label{ehrsimplrecgenfct}
  \Ehr_{ \Delta^\circ } (z) = \frac{ z^{ d+1 } h \left( \frac 1 z \right) }{ (1-z)^{ d+1 } } \, .
\end{equation}
We can compute $\Ehr_{ \Delta^\circ } (z)$ along the same lines as we computed $\Ehr_{ \Delta } (z)$ in part (a):
\[
  \Ehr_{ \Delta^\circ } (z) = \sum_{ t \ge 0 } \# \left( \text{lattice points in } \cone(\Delta^\circ) \text{ at height } t \right) z^t .
\]
The fundamental parallelepiped of $\cone(\Delta^\circ) = \sum_{ \v \text{ vertex of } \Delta } \R_{ > 0 } (\v, 1)$ is
\[
  \widetilde \Q := \sum_{ \v \text{ vertex of } \Delta } (0,1] (\v, 1) \, ,
\]
and
$
  \Ehr_{ \Delta^\circ } (z) = \frac{ \widetilde h(z) }{ (1-z)^{ d+1 } }
$
where
\[
  \widetilde h(z) := \sum_{ t \ge 0 } \# \left( \text{lattice points in } \widetilde \Q \text{ at height } t \right) z^t .
\]
\def\JPicScale{0.7}
\begin{figure}[htb]
\begin{center}
\ifx\JPicScale\undefined\def\JPicScale{1}\fi
\unitlength \JPicScale mm
\begin{picture}(170,50.62)(0,0)
\linethickness{0.3mm}
\put(45,35){\circle*{2.5}}

\linethickness{0.3mm}
\put(5,25){\line(1,0){45}}
\put(50,25){\vector(1,0){0.12}}
\linethickness{0.3mm}
\put(25,0){\line(0,1){50}}
\put(25,50){\vector(0,1){0.12}}
\linethickness{0.3mm}
\put(15,35){\circle*{2.5}}

\linethickness{0.3mm}
\multiput(25,25)(0.24,0.12){83}{\line(1,0){0.24}}
\linethickness{0.3mm}
\multiput(15,35)(0.12,-0.12){83}{\line(1,0){0.12}}
\linethickness{0.3mm}
\multiput(15,35)(1.74,0.87){12}{\multiput(0,0)(0.22,0.11){4}{\line(1,0){0.22}}}
\linethickness{0.3mm}
\multiput(35,45)(1.33,-1.33){8}{\multiput(0,0)(0.11,-0.11){6}{\line(1,0){0.11}}}
\put(53.75,34.38){\makebox(0,0)[cc]{$(\v_2,1)$}}

\put(6.25,34.38){\makebox(0,0)[cc]{$(\v_1,1)$}}

\linethickness{0.3mm}
\put(65,25){\line(1,0){45}}
\put(110,25){\vector(1,0){0.12}}
\linethickness{0.3mm}
\put(85,0){\line(0,1){50}}
\put(85,50){\vector(0,1){0.12}}
\linethickness{0.3mm}
\multiput(65,15)(0.24,0.12){83}{\line(1,0){0.24}}
\linethickness{0.3mm}
\multiput(85,25)(0.12,-0.12){83}{\line(1,0){0.12}}
\linethickness{0.3mm}
\multiput(75,5)(1.74,0.87){12}{\multiput(0,0)(0.22,0.11){4}{\line(1,0){0.22}}}
\linethickness{0.3mm}
\multiput(65,15)(1.33,-1.33){8}{\multiput(0,0)(0.11,-0.11){6}{\line(1,0){0.11}}}
\linethickness{0.3mm}
\put(125,25.62){\line(1,0){45}}
\put(170,25.62){\vector(1,0){0.12}}
\linethickness{0.3mm}
\put(145,0.62){\line(0,1){50}}
\put(145,50.62){\vector(0,1){0.12}}
\linethickness{0.3mm}
\multiput(145,25.62)(1.74,0.87){12}{\multiput(0,0)(0.22,0.11){4}{\line(1,0){0.22}}}
\linethickness{0.3mm}
\multiput(135,35.62)(1.33,-1.33){8}{\multiput(0,0)(0.11,-0.11){6}{\line(1,0){0.11}}}
\linethickness{0.3mm}
\multiput(135,35.62)(0.24,0.12){83}{\line(1,0){0.24}}
\linethickness{0.3mm}
\multiput(155,45.62)(0.12,-0.12){83}{\line(1,0){0.12}}
\put(30,35.62){\makebox(0,0)[cc]{$\Q$}}

\put(78.75,15){\makebox(0,0)[cc]{$-\Q$}}

\put(146.25,20){\makebox(0,0)[cc]{$-\Q+(\v_1,1)+(\v_2,1)$}}

\end{picture}
\end{center}
\caption{An instance of \eqref{fundparreciprocityeq}.}\label{fundparreciprocityfig}
\end{figure}
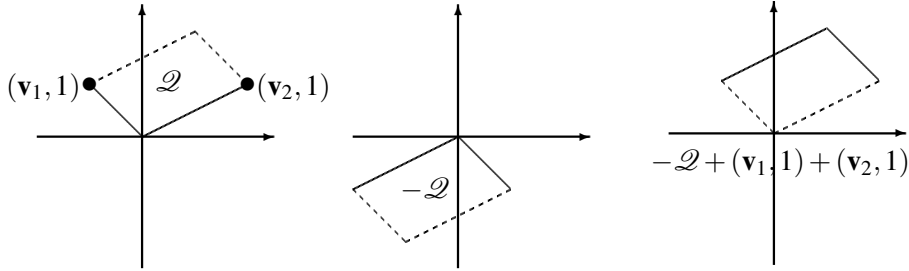

\noindent
Fortunately, the parallelepipeds $\Q$ and $\widetilde \Q$ are geometrically closely related:
\begin{equation}\label{fundparreciprocityeq}
  \widetilde \Q = - \Q + \sum_{ \v \text{ vertex of } \Delta } (\v, 1) \, .
\end{equation}
(Figure \ref{fundparreciprocityfig} shows one instance of this relation.)
This translates into the generating-function relation
\[
  \widetilde h (z) = h \left( \tfrac 1 z \right) z^{ d+1 } 
\]
which proves \eqref{ehrsimplrecgenfct} and thus the second part of Theorem \ref{ehrmacthm} in the simplex case.


\subsection{Lattice Polytopes}\label{ehrhartgeneralsection}

\begin{figure}[htb]
\begin{center}
\includegraphics[totalheight=1.6in]{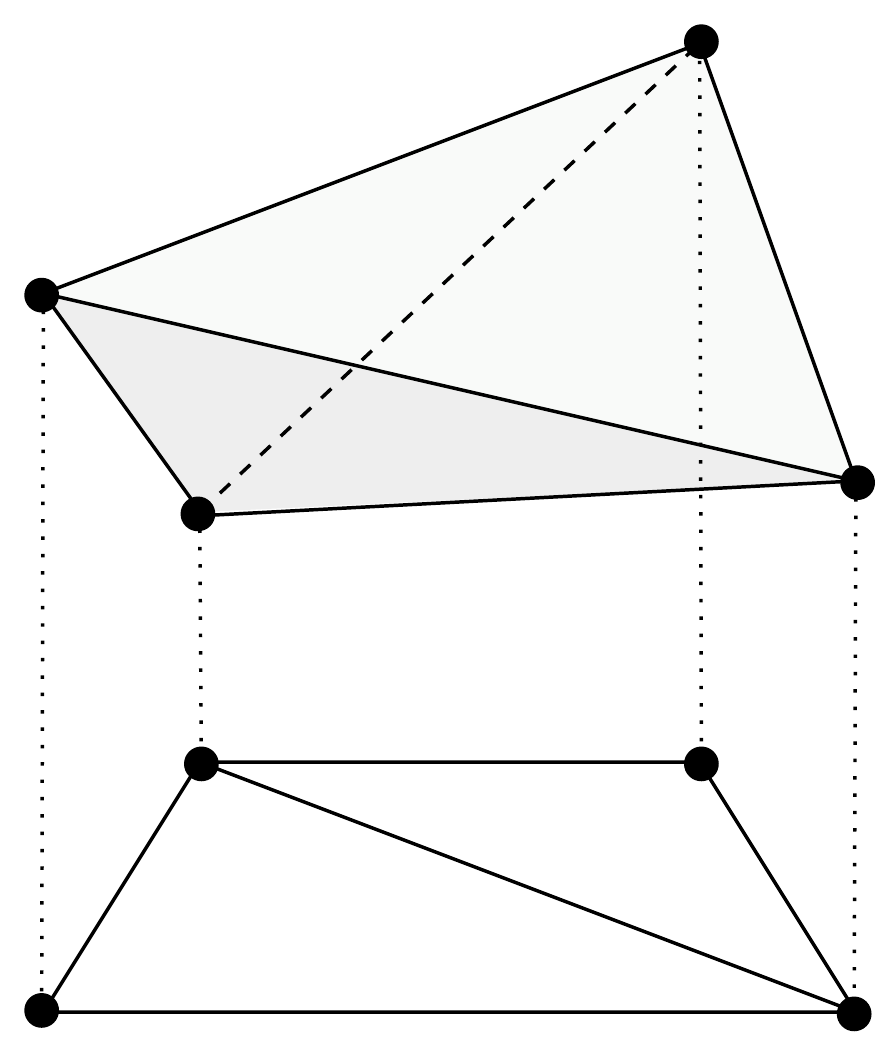}
\end{center}
\caption{A regular triangulation of a quadrilateral.}\label{triangulationfig}
\end{figure}

The general case of Theorem \ref{ehrmacthm} follows from decomposing a general lattice polytope into lattice simplices:
a \emph{triangulation} of a convex $d$-polytope $\P$ is a finite collection $T$ of $d$-simplices with the properties:
\begin{itemize}
\item $\P = \bigcup_{ \Delta \in T } \Delta \, .$
\item For any $\Delta_1, \Delta_2 \in T$, $\Delta_1 \cap \Delta_2$ is a face of both $\Delta_1$ and $\Delta_2$.
\end{itemize}
Here is an algorithm to obtain what's called a \emph{regular triangulation} of $\P = \conv \left\{ \v_1, \v_2, \dots, \v_n \right\} \subset \R^d$:
\begin{enumerate}[(i)]
\item Embed $\P$ into $\R^{d+1}$ as $\P \times \left\{ 0 \right\} .$
\item Randomly choose $r_1, r_2, \dots, r_n \in \R$.
\item Project the lower facets of $\Q := \conv \left\{ (\v_1, r_1), (\v_2, r_2), \dots, (\v_n, r_n) \right\}$ onto $\R^d \times \left\{ 0 \right\} .$
\end{enumerate}
By \emph{lower facets} of $\Q$ we mean those facets that one can see ``from below," i.e., those facets of $\Q$ visible from the point $(\0, -r)$, for some sufficiently large $r$.
Figure \ref{triangulationfig} illustrates the process of obtaining such a triangulation for a quadrilateral.
It's a good exercise to prove that the above algorithm indeed yields a triangulation, for any polytope.\footnote{Strictly speaking, this algorithm yields a triangulation ``only" with probability~1.}
This implies, in particular, that every polytope admits a triangulation (whose simplices have vertices among the vertices of the polytope).

The first part of Theorem \ref{ehrmacthm} follows now immediately, since for a given lattice polytope $\P$ we can write $\ehr_\P(t)$ as a sum/difference of the Ehrhart polynomials of the simplices of a triangulation of $\P$ and their faces, in an inclusion--exclusion way.
To prove the second part of Theorem \ref{ehrmacthm}, we need to work a little harder.
Fix a triangulation of $\P$ and consider the poset $\Phi$ of all faces (including $\emptyset$) of the simplices in this triangulation, ordered by set inclusion.
It will be useful to make $\Phi$ into a lattice, so let's introduce an artificial largest element $\1 \in \Phi$ whose dimension we declare to be $d+1$.
It's a fun (and not entirely trivial) exercise to show that the M\"obius function of $\Phi$ is (assuming that $\G \subset \F$)
\begin{equation}\label{triangulationmobiusthm}
  \mu (\G, \F) = \begin{cases}
    0 & \text{ if } ( \G \subset \partial \P \text{ or } \G = \emptyset ) \text{ and } \F = \1 , \\
    (-1)^{ \dim \F - \dim \G } & \text{ otherwise. }
  \end{cases}
\end{equation}
(Here $\partial \P$ denotes the boundary of $\P$.)
We can now show how the general Ehrhart--Macdonal reciprocity follows from the simplex case. We will use M\"obius inversion \eqref{mobiusinversion} on $\Phi$ for the functions
\[
  f(\F) = \begin{cases}
    \ehr_\F(t) & \text{ if } \F \ne \1, \\
    \ehr_\P(t) & \text{ if } \F = \1,
  \end{cases}
  \qquad \text{ and } \qquad
  g(\F) = \begin{cases}
    \ehr_{\F^\circ}(t) & \text{ if } \F \ne \1, \\
    0 & \text{ if } \F = \1.
  \end{cases}
\]
Because every point in $\P$ is in the interior of a unique face,\footnote{
Here we mean \emph{relative} interior; in particular, $\F^\circ = \F$ if $\F$ is a vertex.
}
\[
  f(\1)
  = \ehr_\P(t)
  = \sum_{ \F \in \Phi \setminus \{ \1 \} } \ehr_{ \F^\circ } (t)
  = \sum_{ \F \in \Phi } g(\F) \, .
\]
By \eqref{mobiusinversion} and \eqref{triangulationmobiusthm},
\[
  g(\1)
  = 0
  = \sum_{ \F \in \Phi } \mu(\F, \1) \, f(\F)
  = \ehr_\P(t) + \sum_{ {\F \in \Phi \setminus \{ \1 \} } \atop {\F \not\subset \partial \P} } (-1)^{ d+1-\dim \F } \ehr_\F(t) \, ,
\]
that is,
\[
  \ehr_\P(t)
  = (-1)^d \sum_{ {\F \in \Phi \setminus \{ \1 \} } \atop {\F \not\subset \partial \P} } (-1)^{ \dim \F } \ehr_\F(t) \, .
\]
Now we evaluate these polynomials at negative integers and use Ehrhart--Macdonald reciprocity for the simplices $\F \in \Phi \setminus \{ \1 \}$:
\[
  \ehr_\P(-t)
  = (-1)^d \sum_{ {\F \in \Phi \setminus \{ \1 \} } \atop {\F \not\subset \partial \P} } (-1)^{ \dim \F } \ehr_\F(-t)
  = (-1)^d \sum_{ {\F \in \Phi \setminus \{ \1 \} } \atop {\F \not\subset \partial \P} } \ehr_{\F^\circ}(t)
  = (-1)^d \ehr_{ \P^\circ } (t) \, ,
\]
and this concludes our proof of Theorem~\ref{ehrmacthm}.

Ehrhart theory is not limited to lattice polytopes; we can relax the integrality condition on the coordinates of the vertices of $\P$ to the rational case. Then $\ehr_\P(t)$ becomes a \emph{quasipolynomial}, i.e., a function of the form
\[
  c_{n}(t) \ t^{n} + \dots + c_{1}(t) \ t + c_{0}(t) \, ,
\]
where $c_{0}, c_1, \dots , c_{n}$ are periodic functions in $t$.
Ehrhart--Macdonald reciprocity carries over verbatim to the rational case.
Further yet, very recent results \cite{baldoniberlinekoeppevergnerealdilation,barvinokfixfirstcoeffs,linkerational} extended Ehrhart (quasi-)polynomials by allowing \emph{rational} or \emph{real} dilation factors when counting lattice points in rational polytopes.

We finish this section by mentioning that there are alternative ways of proving Theorem \ref{ehrmacthm}, see, e.g., \cite[Chapter 4]{ccd} and \cite{samehrhartpolynomial}; our proof followed Ehrhart's original lines \cite{ehrhartpolynomial} (Section \ref{ehrhartsimplexsection}) and \cite[Chapter 4]{stanleyec1} (Section \ref{ehrhartgeneralsection}).
For (much) more about triangulations, we recommend \cite{deloerarambausantos}; for more about Ehrhart polynomials, see \cite{ccd}, \cite{hibi}, and~\cite[Chapter 4]{stanleyec1}.


\section{A Polyhedral View at Graph Colorings and Acyclic Orientations}\label{coloringiopsection}

Our next step is to interpret graph coloring geometrically, with the goal of deriving Theorem \ref{coloringreciprocitythm}.
After having meditated about lattice point in polytopes for a while now, it is a short step to view a coloring $\x \in [t]^V$ of a graph $G=(V,E)$ as an integer point in the cube $[1,t]^V$ or, more conveniently, an \emph{interior} lattice point in the $(t+1)$-dilate of the unit cube $[0,1]^V$.
This $t$-coloring $\x$ is proper if it misses the hyperplane arrangement
\[
  \H_G := \left\{ x_i = x_j : \, ij \in E \right\} ,
\]
the \emph{graphical arrangement} corresponding to $G$.
Thus each proper $t$-coloring corresponds to a lattice point in 
\begin{equation}\label{coloringiopeq}
  \left( (t+1) \P^\circ \setminus \bigcup \H_G \right) \cap \Z^V ,
\end{equation}
where $\P = [0,1]^V$ is the unit cube in $\R^V$ (see the left-hand side of Figure \ref{F:gcol2} for an example where $G = K_2$, the graph with exactly two adjacent nodes).

\def\JPicScale{0.5}
\begin{figure}[htb]
\begin{center}
\ifx\JPicScale\undefined\def\JPicScale{1}\fi
\unitlength \JPicScale mm
\begin{picture}(235,94)(0,0)
\linethickness{0.3mm}
\put(30,20){\circle*{2.5}}

\linethickness{0.3mm}
\put(10,0){\line(0,1){90}}
\put(10,90){\vector(0,1){0.12}}
\linethickness{0.3mm}
\put(0,10){\line(1,0){90}}
\put(90,10){\vector(1,0){0.12}}
\linethickness{0.3mm}
\multiput(0,0)(1.42,1.42){64}{\multiput(0,0)(0.12,0.12){6}{\line(1,0){0.12}}}
\linethickness{0.3mm}
\multiput(10,80)(1.97,0){36}{\line(1,0){0.99}}
\linethickness{0.3mm}
\multiput(80,10)(0,1.97){36}{\line(0,1){0.99}}
\put(101,93){\makebox(0,0)[cc]{$x_1=x_2$}}

\put(95,90){\makebox(0,0)[cc]{}}

\linethickness{0.3mm}
\put(50,20){\circle*{2.5}}

\linethickness{0.3mm}
\put(40,20){\circle*{2.5}}

\linethickness{0.3mm}
\put(40,30){\circle*{2.5}}

\linethickness{0.3mm}
\put(60,30){\circle*{2.5}}

\linethickness{0.3mm}
\put(50,30){\circle*{2.5}}

\linethickness{0.3mm}
\put(50,40){\circle*{2.5}}

\linethickness{0.3mm}
\put(70,40){\circle*{2.5}}

\linethickness{0.3mm}
\put(60,40){\circle*{2.5}}

\linethickness{0.3mm}
\put(70,20){\circle*{2.5}}

\linethickness{0.3mm}
\put(60,20){\circle*{2.5}}

\linethickness{0.3mm}
\put(70,30){\circle*{2.5}}

\linethickness{0.3mm}
\put(70,50){\circle*{2.5}}

\linethickness{0.3mm}
\put(60,50){\circle*{2.5}}

\linethickness{0.3mm}
\put(70,60){\circle*{2.5}}

\put(95,9.38){\makebox(0,0)[cc]{$x_1$}}

\put(10,94){\makebox(0,0)[cc]{$x_2$}}

\put(15,90){\makebox(0,0)[cc]{}}

\linethickness{0.3mm}
\put(20,70){\circle*{2.5}}

\linethickness{0.3mm}
\put(40,70){\circle*{2.5}}

\linethickness{0.3mm}
\put(30,70){\circle*{2.5}}

\linethickness{0.3mm}
\put(60,70){\circle*{2.5}}

\linethickness{0.3mm}
\put(50,70){\circle*{2.5}}

\linethickness{0.3mm}
\put(20,60){\circle*{2.5}}

\linethickness{0.3mm}
\put(40,60){\circle*{2.5}}

\linethickness{0.3mm}
\put(30,60){\circle*{2.5}}

\linethickness{0.3mm}
\put(50,60){\circle*{2.5}}

\linethickness{0.3mm}
\put(20,50){\circle*{2.5}}

\linethickness{0.3mm}
\put(40,50){\circle*{2.5}}

\linethickness{0.3mm}
\put(30,50){\circle*{2.5}}

\linethickness{0.3mm}
\put(20,40){\circle*{2.5}}

\linethickness{0.3mm}
\put(30,40){\circle*{2.5}}

\linethickness{0.3mm}
\put(20,30){\circle*{2.5}}

\put(80,5){\makebox(0,0)[cc]{}}

\put(80,4){\makebox(0,0)[cc]{$t+1$}}

\put(1,80){\makebox(0,0)[cc]{$t+1$}}

\linethickness{0.3mm}
\put(80,8.75){\line(0,1){2.5}}
\linethickness{0.3mm}
\put(8.75,80){\line(1,0){2.5}}
\linethickness{0.3mm}
\put(160,10){\circle*{2.5}}

\linethickness{0.3mm}
\put(180,10){\circle*{2.5}}

\linethickness{0.3mm}
\put(170,10){\circle*{2.5}}

\linethickness{0.3mm}
\put(170,20){\circle*{2.5}}

\linethickness{0.3mm}
\put(190,20){\circle*{2.5}}

\linethickness{0.3mm}
\put(180,20){\circle*{2.5}}

\linethickness{0.3mm}
\put(180,30){\circle*{2.5}}

\linethickness{0.3mm}
\put(200,30){\circle*{2.5}}

\linethickness{0.3mm}
\put(190,30){\circle*{2.5}}

\linethickness{0.3mm}
\put(200,10){\circle*{2.5}}

\linethickness{0.3mm}
\put(190,10){\circle*{2.5}}

\linethickness{0.3mm}
\put(200,20){\circle*{2.5}}

\linethickness{0.3mm}
\put(200,40){\circle*{2.5}}

\linethickness{0.3mm}
\put(190,40){\circle*{2.5}}

\linethickness{0.3mm}
\put(200,50){\circle*{2.5}}

\linethickness{0.3mm}
\put(150,60){\circle*{2.5}}

\linethickness{0.3mm}
\put(170,60){\circle*{2.5}}

\linethickness{0.3mm}
\put(160,60){\circle*{2.5}}

\linethickness{0.3mm}
\put(190,60){\circle*{2.5}}

\linethickness{0.3mm}
\put(180,60){\circle*{2.5}}

\linethickness{0.3mm}
\put(150,50){\circle*{2.5}}

\linethickness{0.3mm}
\put(170,50){\circle*{2.5}}

\linethickness{0.3mm}
\put(160,50){\circle*{2.5}}

\linethickness{0.3mm}
\put(180,50){\circle*{2.5}}

\linethickness{0.3mm}
\put(150,40){\circle*{2.5}}

\linethickness{0.3mm}
\put(170,40){\circle*{2.5}}

\linethickness{0.3mm}
\put(160,40){\circle*{2.5}}

\linethickness{0.3mm}
\put(150,30){\circle*{2.5}}

\linethickness{0.3mm}
\put(160,30){\circle*{2.5}}

\linethickness{0.3mm}
\put(150,20){\circle*{2.5}}

\linethickness{0.3mm}
\put(190,50){\circle*{2.5}}

\linethickness{0.3mm}
\put(180,40){\circle*{2.5}}

\linethickness{0.3mm}
\put(170,30){\circle*{2.5}}

\linethickness{0.3mm}
\put(160,20){\circle*{2.5}}

\linethickness{0.3mm}
\put(150,10){\circle*{2.5}}

\linethickness{0.3mm}
\put(200,60){\circle*{2.5}}

\linethickness{0.3mm}
\put(150,60){\line(1,0){51}}
\linethickness{0.3mm}
\put(200,10){\line(0,1){50}}
\linethickness{0.3mm}
\multiput(150,10.38)(0.12,0.12){417}{\line(1,0){0.12}}
\linethickness{0.3mm}
\multiput(150,9.75)(0.12,0.12){417}{\line(1,0){0.12}}
\put(200,4){\makebox(0,0)[cc]{$t-1$}}

\put(140,60){\makebox(0,0)[cc]{$t-1$}}

\linethickness{0.3mm}
\put(150,0){\line(0,1){90}}
\put(150,90){\vector(0,1){0.12}}
\linethickness{0.3mm}
\put(140,10){\line(1,0){90}}
\put(230,10){\vector(1,0){0.12}}
\put(235,9.38){\makebox(0,0)[cc]{$x_1$}}

\put(150,94){\makebox(0,0)[cc]{$x_2$}}

\end{picture}
\end{center}
\caption{The integer points $t$-color the graph $K_2$ (with $t=6$) and the reciprocal picture.} \label{F:gcol2}
\end{figure}
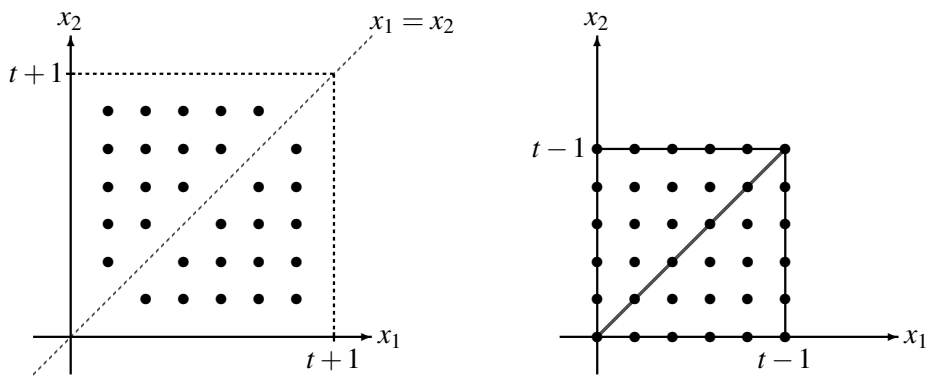

Viewed like this, counting proper $t$-colorings is quite reminiscent of Ehrhart theory, safe for the graphic arrangement whose hyperplanes contain the non-proper colorings.
At any rate, $\P^\circ \setminus \bigcup \H_G$ is a union of open polytopes, say
\[
  \P^\circ \setminus \bigcup \H_G = \Q_1^\circ \cup \Q_2^\circ \cup \dots \cup \Q_n^\circ \, ,
\]
and so we can indeed express the chromatic polynomial in Ehrhartian terms:
\[
  \chrom_G(t) = \sum_{ j=1 }^n \ehr_{ \Q_j^\circ } (t+1) \, .
\]
The reciprocal counting function is therefore, by Ehrhart--Macdonald reciprocity (Theorem \ref{ehrmacthm}),
\[
  \chrom_G(-t) = (-1)^{ |V| }  \sum_{ j=1 }^n \ehr_{ \Q_j } (t-1)
\]
since all $\Q_j$'s have the same dimension $|V|$.
(On the right in Figure \ref{F:gcol2} is an illustration of this count for $G = K_2$.)
The right-hand side counts lattice points in the (closed) cube $(t-1) \P$ with multiplicity: each lattice point $\x$ gets weighted by the number of $\Q_j$'s containing it; geometrically this is the number of closed regions of $\H_G$ containing $\x$.
The last ingredient for our proof of Theorem \ref{coloringreciprocitythm} is the following simple but crucial observation, illustrated in Figure~\ref{acyclicorientfig}.

\def\JPicScale{0.6}
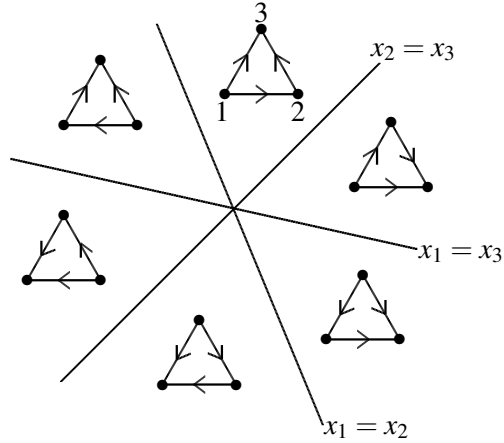
\begin{figure}[htb]
\begin{center}
\ifx\JPicScale\undefined\def\JPicScale{1}\fi
\unitlength \JPicScale mm
\begin{picture}(100,92.5)(0,0)
\linethickness{0.3mm}
\put(41.87,24.37){\circle*{2.5}}

\linethickness{0.3mm}
\put(50,10){\circle*{2.5}}

\linethickness{0.3mm}
\put(33.75,10){\circle*{2.5}}

\linethickness{0.3mm}
\multiput(0,60)(0.54,-0.12){167}{\line(1,0){0.54}}
\linethickness{0.3mm}
\multiput(32.5,90)(0.12,-0.29){307}{\line(0,-1){0.29}}
\linethickness{0.3mm}
\multiput(11.25,10.62)(0.12,0.12){589}{\line(0,1){0.12}}
\linethickness{0.3mm}
\multiput(41.88,24.38)(0.12,-0.21){68}{\line(0,-1){0.21}}
\linethickness{0.3mm}
\put(33.75,10){\line(1,0){16.25}}
\linethickness{0.3mm}
\multiput(33.75,10)(0.12,0.21){68}{\line(0,1){0.21}}
\put(79,0){\makebox(0,0)[cc]{$x_1=x_2$}}

\put(100,39){\makebox(0,0)[cc]{$x_1=x_3$}}

\put(89.38,83.75){\makebox(0,0)[cc]{$x_2=x_3$}}

\linethickness{0.3mm}
\multiput(40.62,10)(0.16,-0.12){16}{\line(1,0){0.16}}
\linethickness{0.3mm}
\multiput(40.62,10)(0.16,0.12){16}{\line(1,0){0.16}}
\linethickness{0.3mm}
\put(36.88,15.62){\line(0,1){3.13}}
\linethickness{0.3mm}
\multiput(36.88,15.62)(0.31,0.13){10}{\line(1,0){0.31}}
\linethickness{0.3mm}
\put(46.88,15.62){\line(0,1){3.13}}
\linethickness{0.3mm}
\multiput(43.75,16.88)(0.31,-0.13){10}{\line(1,0){0.31}}
\linethickness{0.3mm}
\put(78.13,34.37){\circle*{2.5}}

\linethickness{0.3mm}
\put(86.25,20){\circle*{2.5}}

\linethickness{0.3mm}
\put(70,20){\circle*{2.5}}

\linethickness{0.3mm}
\multiput(78.12,34.38)(0.12,-0.21){68}{\line(0,-1){0.21}}
\linethickness{0.3mm}
\put(70,20){\line(1,0){16.25}}
\linethickness{0.3mm}
\multiput(70,20)(0.12,0.21){68}{\line(0,1){0.21}}
\linethickness{0.3mm}
\put(73.12,25.62){\line(0,1){3.13}}
\linethickness{0.3mm}
\multiput(73.12,25.62)(0.31,0.13){10}{\line(1,0){0.31}}
\linethickness{0.3mm}
\put(83.12,25.62){\line(0,1){3.13}}
\linethickness{0.3mm}
\multiput(80,26.88)(0.31,-0.13){10}{\line(1,0){0.31}}
\linethickness{0.3mm}
\multiput(76.88,21.88)(0.16,-0.12){16}{\line(1,0){0.16}}
\linethickness{0.3mm}
\multiput(76.88,18.12)(0.16,0.12){16}{\line(1,0){0.16}}
\linethickness{0.3mm}
\put(11.87,47.5){\circle*{2.5}}

\linethickness{0.3mm}
\put(20,33.13){\circle*{2.5}}

\linethickness{0.3mm}
\put(3.75,33.13){\circle*{2.5}}

\linethickness{0.3mm}
\multiput(11.88,47.5)(0.12,-0.21){68}{\line(0,-1){0.21}}
\linethickness{0.3mm}
\put(3.75,33.12){\line(1,0){16.25}}
\linethickness{0.3mm}
\multiput(3.75,33.12)(0.12,0.21){68}{\line(0,1){0.21}}
\linethickness{0.3mm}
\multiput(10.62,33.12)(0.16,-0.12){16}{\line(1,0){0.16}}
\linethickness{0.3mm}
\multiput(10.62,33.12)(0.16,0.12){16}{\line(1,0){0.16}}
\linethickness{0.3mm}
\put(6.88,38.75){\line(0,1){3.13}}
\linethickness{0.3mm}
\multiput(6.88,38.75)(0.31,0.12){10}{\line(1,0){0.31}}
\linethickness{0.3mm}
\multiput(15.62,41.25)(0.16,-0.12){16}{\line(1,0){0.16}}
\linethickness{0.3mm}
\put(15.62,38.12){\line(0,1){3.13}}
\linethickness{0.3mm}
\put(84.37,68.13){\circle*{2.5}}

\linethickness{0.3mm}
\put(92.5,53.75){\circle*{2.5}}

\linethickness{0.3mm}
\put(76.25,53.75){\circle*{2.5}}

\linethickness{0.3mm}
\multiput(84.38,68.12)(0.12,-0.21){68}{\line(0,-1){0.21}}
\linethickness{0.3mm}
\put(76.25,53.75){\line(1,0){16.25}}
\linethickness{0.3mm}
\multiput(76.25,53.75)(0.12,0.21){68}{\line(0,1){0.21}}
\linethickness{0.3mm}
\put(89.38,59.38){\line(0,1){3.12}}
\linethickness{0.3mm}
\multiput(86.25,60.62)(0.31,-0.12){10}{\line(1,0){0.31}}
\linethickness{0.3mm}
\multiput(83.12,55.62)(0.16,-0.12){16}{\line(1,0){0.16}}
\linethickness{0.3mm}
\multiput(83.12,51.88)(0.16,0.12){16}{\line(1,0){0.16}}
\linethickness{0.3mm}
\put(81.25,58.75){\line(0,1){3.75}}
\linethickness{0.3mm}
\multiput(78.12,60.62)(0.2,0.12){16}{\line(1,0){0.2}}
\linethickness{0.3mm}
\put(55.63,88.75){\circle*{2.5}}

\linethickness{0.3mm}
\put(63.75,74.37){\circle*{2.5}}

\linethickness{0.3mm}
\put(47.5,74.37){\circle*{2.5}}

\linethickness{0.3mm}
\multiput(55.62,88.75)(0.12,-0.21){68}{\line(0,-1){0.21}}
\linethickness{0.3mm}
\put(47.5,74.38){\line(1,0){16.25}}
\linethickness{0.3mm}
\multiput(47.5,74.38)(0.12,0.21){68}{\line(0,1){0.21}}
\linethickness{0.3mm}
\multiput(54.38,76.25)(0.16,-0.12){16}{\line(1,0){0.16}}
\linethickness{0.3mm}
\multiput(54.38,72.5)(0.16,0.12){16}{\line(1,0){0.16}}
\linethickness{0.3mm}
\put(52.5,79.38){\line(0,1){3.74}}
\linethickness{0.3mm}
\multiput(49.38,81.25)(0.19,0.12){16}{\line(1,0){0.19}}
\linethickness{0.3mm}
\put(58.75,79.38){\line(0,1){3.74}}
\linethickness{0.3mm}
\multiput(58.75,83.12)(0.2,-0.12){16}{\line(1,0){0.2}}
\linethickness{0.3mm}
\put(20,81.87){\circle*{2.5}}

\linethickness{0.3mm}
\put(28.13,67.5){\circle*{2.5}}

\linethickness{0.3mm}
\put(11.87,67.5){\circle*{2.5}}

\linethickness{0.3mm}
\multiput(20,81.88)(0.12,-0.21){68}{\line(0,-1){0.21}}
\linethickness{0.3mm}
\put(11.88,67.5){\line(1,0){16.24}}
\linethickness{0.3mm}
\multiput(11.88,67.5)(0.12,0.21){68}{\line(0,1){0.21}}
\linethickness{0.3mm}
\put(16.88,72.5){\line(0,1){3.75}}
\linethickness{0.3mm}
\multiput(13.75,74.38)(0.2,0.12){16}{\line(1,0){0.2}}
\linethickness{0.3mm}
\put(23.12,72.5){\line(0,1){3.75}}
\linethickness{0.3mm}
\multiput(23.12,76.25)(0.2,-0.12){16}{\line(1,0){0.2}}
\linethickness{0.3mm}
\multiput(18.75,67.5)(0.16,0.12){16}{\line(1,0){0.16}}
\linethickness{0.3mm}
\multiput(18.75,67.5)(0.16,-0.12){16}{\line(1,0){0.16}}
\put(46.88,70.62){\makebox(0,0)[cc]{$1$}}

\put(63.75,70.62){\makebox(0,0)[cc]{$2$}}

\put(55.62,92.5){\makebox(0,0)[cc]{$3$}}

\end{picture}
\end{center}
\caption{The regions of $\H_{ K_3 }$ (projected to the plane $x_1+x_2+x_3 = 0$) and their corresponding acyclic orientations.} \label{acyclicorientfig}
\end{figure}

\begin{lemma}[Greene \cite{greeneacyclic,greenezaslavsky}]\label{greenelemma}
The regions of $\H_G$ are in one-to-one correspondence with the acyclic orientations of~$G$.
\end{lemma}

Theorem \ref{coloringreciprocitythm} follows now by (re-)interpreting the lattice points in $(t-1) \P$ as $t$-colorings and interpreting their multiplicities in terms of compatible acyclic orientations.


\section{Inside-out Polytopes}\label{iopsection}

The above proof of Theorem \ref{coloringreciprocitythm} appeared in \cite{iop}; we take a short detour to illustrate how other reciprocity theorems follow from this work.
The scenery of our proof consisted of a (rational) polytope $\P$, a (rational) hyperplane arrangement $\H$, and the two counting functions\footnote{The shift from dilating polytopes to shrinking the lattice is purely cosmetic, as there may be hyperplanes in $\H$ that do not contain the origin.}
\[
  I_{ \P, \H } (t) := \# \left( \left( \P \setminus \bigcup \H \right) \cap \frac 1 t \Z^d \right)
\qquad \text{ and } \qquad
  O_{ \P, \H } (t) := \sum_{ \x \in \frac 1 t \Z^d } \mult_{ \P, \H } (\x)
\]
where
\[
  \mult_{ \P, \H } (\x) := \begin{cases}
    \text{number of closed regions of $(\P, \H)$ that contain $\x$} & \text{ if } \x \in \P,
    \\
    0 & \text{ if } \x \notin \P .
  \end{cases}
\]
The pair $(\P, \H)$ goes by the name \emph{inside-out polytope} (we think of the hyperplanes in $\H$ as acting as additional boundary of the polytope $\P$ ``turned inside out"), and our above application of Ehrhart--Macdonald reciprocity (Theorem \ref{ehrmacthm}) shows that the two inside-out polytope counting functions are reciprocal quasipolynomials~\cite{iop}:
\begin{equation}\label{ioprec}
  I_{ \P^\circ, \H } (-t) = (-1)^{\dim \P} O_{ \P, \H } (t) \, .
\end{equation}
Looking back once more at our above proof of Theorem \ref{coloringreciprocitythm} illustrates the two central ingredients we need in order to apply \eqref{ioprec} to a specific combinatorial situation: first, we need to be able to interpret the underlying objects that we are counting as lattice points in $t \P^\circ \setminus \H$ (or some close variant); once we have this interpretation, we can apply \eqref{ioprec}, in other words, we are guaranteed a reciprocity theorem in the world of polyhedral geometry. The ``big question" is whether we can return into the world of the original combinatorial situation, in other words, if we can interpret the multiplicities appearing in $O_{ \P, \H } (t)$ in that world. In the graph-coloring case, this last step was made possible by Lemma \ref{greenelemma}; the ``big question" we just mentioned thus reduces essentially to finding an analogous result in the given combinatorial situation.

Fortunately, there is a number of combinatorial constructs in which the inside-out polytope approach resulted in (novel) reciprocity theorems \cite{nowhereharmonic,magiclatin,iop,nnz,breuerdall,breuersanyal,zaslavskybiasedgraphs7}.



\section{A Polyhedral View at P-partitions}\label{ppartsection}

In the previous two sections, we arranged Ehrhart (quasi-)polynomials with hyperplanes, in the sense that we enumerated lattice points in polyhedra but excluded lattice points on certain hyperplanes.
We will now exhibit a second mix of Ehrhart theory and hyperplane arrangements: we will use hyperplanes to triangulate polyhedra whose lattice points we want to enumerate.

Suppose $\Pi = \left\{ a_1, a_2, \dots, a_d \right\}$ is a poset. For technical reasons which will become clear soon, we assume that the indices of the $a_j$'s respect the order of $\Pi$ in the sense that we have $j \le k$ if $a_j \preceq a_k$. For example, we need to re-lable Figure \ref{introposetfig} in such a way that $b = a_3$ (because $b$ is the maximal element in this poset).
With this convention and in sync with Section \ref{ppartintrosection}, we define the set of all $\Pi$-partitions as
\[
  \Kappa(\Pi) := \left\{ \x \in \Z_{ \ge 0 }^d : \, x_j \ge x_k \ \text{ if } \ a_j \preceq a_k \right\} .
\]
A \emph{linear extension} of $\Pi$ is a chain $\Gamma$ on $\left\{ a_1, a_2, \dots, a_d \right\}$ that preserves any relation of $\Pi$.
The relations in $\Gamma$ are uniquely determined by a permutation $\sigma \in S_d$, namely the one that orders the chain:
\[
  a_{ \sigma(1) } \prec a_{ \sigma(2) } \prec \dots \prec a_{ \sigma(d) } \, ;
\]
we will call this chain $\Gamma_\sigma$.
Not every permutation $\sigma \in S_d$ will give rise to a linear extension $\Gamma_\sigma$ of $\Pi$, but only those $\sigma$ that respect the order of $\Pi$, i.e., 
\begin{equation}\label{partitionrespectingorder}
  a_j \preceq a_k \ \text{ in } \ \Pi
  \qquad \Longrightarrow \qquad
  a_j \preceq a_k \ \text{ in } \ \Gamma_\sigma \, .
\end{equation}
\def\JPicScale{0.6}
\begin{figure}[htb]
\begin{center}
\input{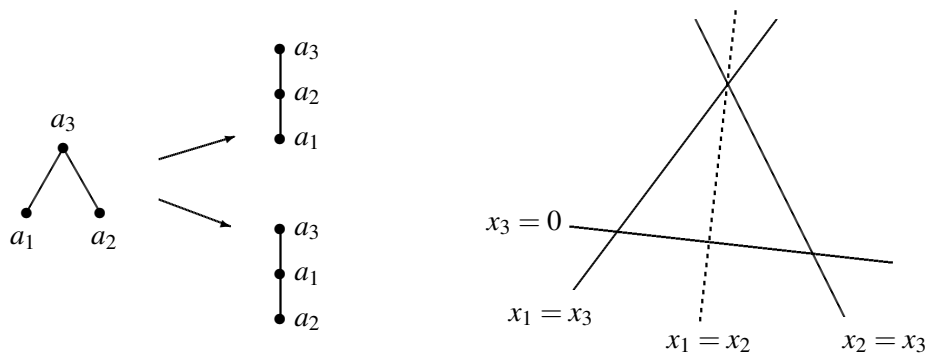}
\end{center}
\caption{The two linear extensions of $\Lambda$, with accompanying cones.}\label{introposetlinextfig}
\end{figure}
For example, the poset $\Lambda$ in Figure \ref{introposetfig} has two linear extensions $\Gamma_\sigma$, for $\sigma = [123]$ and $[213]$ (written in one-line notation), which are pictured on the left in Figure~\ref{introposetlinextfig}.
We can see in this example that
$
  \Kappa(\Lambda)
  = \Kappa \left( \Gamma_{ [123] } \right) \ \cup \ \Kappa \left( \Gamma_{ [213] } \right)
$; more generally, for any poset $\Pi$, we have
\begin{equation}\label{pparttriangulation}
  \Kappa(\Pi) \ = \ \bigcup_\sigma \ \Kappa\left( \Gamma_\sigma \right) ,
\end{equation}
where the union is taken over all $\sigma \in S_d$ that satisfy \eqref{partitionrespectingorder}.
It is natural to think of the elements of $\Kappa(\Pi)$ as lattice points in the cone
\[
  \left\{ \x \in \R_{ \ge 0 }^d : \, x_j \ge x_k \ \text{ if } \ a_j \preceq a_k \right\} ,
\]
and \eqref{pparttriangulation} gives a triangulation of this cone.
On the right in Figure \ref{introposetlinextfig}, we can see how this triangulation looks for the cone behind $\Lambda$ (rather, a two-dimensional slice of this three-dimensional cone).
In fact, we can say more: first, this triangulation is \emph{unimodular}, i.e., each cone represented on the right-hand side of \eqref{pparttriangulation} has generators that span the integer lattice.
Second, we can write \eqref{pparttriangulation} as a \emph{disjoint} union by making use of the \emph{descent set} of a permutation $\sigma \in S_d$, defined as
\[
  \Des \sigma := \left\{ j \in [d-1] : \, \sigma(j) > \sigma(j+1) \right\} .
\]
In our running example $\Pi = \Lambda$, $\Des[123] = \emptyset$ and $\Des[213] = \left\{ 1 \right\}$.
So by writing
\[
  \widetilde\Kappa \left( \Gamma_\sigma \right) := \left\{ \x \in \Z_{ \ge 0 }^d :
     \begin{array}{ll}
       x_{ \sigma(1) } \ge x_{ \sigma(2) } \ge \dots \ge x_{ \sigma(d) } , \\
       x_{ \sigma(j) } > x_{ \sigma(j+1) } \text{ if } j \in \Des \sigma
     \end{array}
  \right\} ,
\]
we obtain
\begin{equation}\label{disjointppartunion}
  \Kappa(\Pi) \ = \ \bigcup_\sigma \ \widetilde\Kappa \left( \Gamma_\sigma \right) ,
\end{equation}
where this now \emph{disjoint} union is taken over all $\sigma \in S_d$ that satisfy \eqref{partitionrespectingorder}.
The generating functions of $\widetilde\Kappa \left( \Gamma_\sigma \right)$ can be computed from first principles with the help of the \emph{major index}
\[
  \maj \sigma := \sum_{ j \in \Des \sigma } j \, ;
\]
it's a fun exercise to show that
\[
  \sum_{ \x \in \widetilde\Kappa \left( \Gamma_\sigma \right) } z^{ x_1 + x_2 + \dots + x_d }
  = \frac{ z^{ \maj \sigma } }{ (1-z) (1-z^2) \cdots (1-z^d) } \, ,
\]
and together with \eqref{disjointppartunion}, this implies:

\begin{lemma}[Stanley \cite{stanleythesis}]\label{ppartlemma}
Let $\Pi$ be a poset on $d$ elements. Then
\[
  P_\Pi(z) = \frac{ \sum_\sigma z^{ \maj \sigma } }{ (1-z) (1-z^2) \cdots (1-z^d) } \, ,
\]
where the sum is taken over all $\sigma \in S_d$ that satisfy \eqref{partitionrespectingorder}.
\end{lemma}

For the analogous lemma for \emph{strict} $\Pi$-partitions, we consider the \emph{ascent set} of a permutation $\sigma \in S_d$,
\[
  \Asc \sigma := \left\{ j \in [d-1] : \, \sigma(j) < \sigma(j+1) \right\} ,
\qquad \text{ and } \qquad
  \amaj \sigma := \sum_{ j \in \Asc \sigma } j \, .
\]
Then
\begin{equation}\label{ppartitionamajdecompex}
  P_\Pi^\circ(z) = \frac{ \sum_\sigma z^{ \amaj \sigma } }{ (1-z) (1-z^2) \cdots (1-z^d) } \, ,
\end{equation}
where once more the sum is taken over all $\sigma \in S_d$ that satisfy \eqref{partitionrespectingorder}.
Theorem \ref{pipartitionreciprocitythm} follows now essentially from the fact that descents and ascents of a permutation are complementary, and so
\begin{equation}\label{desascsum}
  \maj \sigma + \amaj \sigma = \sum_{ j=1 }^{ d-1 } j = \tbinom d 2 \, .
\end{equation}
By Lemma \ref{ppartlemma},
\begin{align*}
  P_\Pi( \tfrac 1 z )
  &= \frac{ \sum_\sigma z^{ -\maj \sigma } }{ (1-z^{ -1 } ) (1-z^{ -2 } ) \cdots (1-z^{ -d } ) }
  = (-1)^d \frac{ z^{ 1 + 2 + \dots + d } \sum_\sigma z^{ - \maj \sigma } }{ (1-z) (1-z^2) \cdots (1-z^d) } \\
  &= (-1)^d \frac{ z^d \sum_\sigma z^{ \binom d 2 - \maj \sigma } }{ (1-z) (1-z^2) \cdots (1-z^d) }
  \stackrel{ \eqref{desascsum} }{=} (-z)^d \frac{ \sum_\sigma z^{ \amaj \sigma } }{ (1-z) (1-z^2) \cdots (1-z^d) }
  \stackrel{ \eqref{ppartitionamajdecompex} }{=} (-z)^{ |\Pi| } P_\Pi^\circ(z) \, ,
\end{align*}
where each sum is taken over all $\sigma \in S_d$ that satisfy \eqref{partitionrespectingorder}.

We close this secion by remarking that Stanley's original approach to $P$-partitions \cite{stanleythesis} is less geometric than our treatment, though one can easily interpret his work along these lines.
The recent papers \cite{mahonianpartition,beckgesselleesavage} used similar discrete-geometric approaches to (number-theoretic) partition identities, where again descent statistics play a role.


\section{Open Problems}

We finish our tour by mentioning a general open problem about all polynomials that appeared as counting functions in this paper, namely the question of classification: give conditions on $a_0, a_1, \dots, a_d$ that allow us to detect whether or not a given polynomial $a_d \, t^d + a_{ d-1 } \, t^{ d-1 } + \dots + a_0$ is a face-number, characteristic, Ehrhart, or chromatic polynomial.
In general, this is a much-too-big research program; for example, the classification problem for Ehrhart polynomials is open already in dimension three.
On the other hand, there has been some exciting recent progress; see, e.g., \cite{haasenillpayne,stapledondelta,stapledonadditive}.
For numerous more open problems about the various combinatorial objects we discussed here, we refer to the books \cite{ccd,deloerarambausantos,stanleyec1,ziegler}.


\bibliographystyle{amsplain}
\bibliography{bib}

\def\cprime{$'$} \def\cprime{$'$}
\providecommand{\bysame}{\leavevmode\hbox to3em{\hrulefill}\thinspace}
\providecommand{\MR}{\relax\ifhmode\unskip\space\fi MR }
\providecommand{\MRhref}[2]{%
  \href{http://www.ams.org/mathscinet-getitem?mr=#1}{#2}
}
\providecommand{\href}[2]{#2}
\begin{thebibliography}{10}

\bibitem{baldoniberlinekoeppevergnerealdilation}
Velleda Baldoni, Nicole Berline, Matthias K{\"o}ppe, and Mich{\`e}le Vergne,
  \emph{Intermediate sums on polyhedra: Computation and real {E}hrhart theory},
  Preprint ({\tt arXiv:1011.6002v1}), 2010.

\bibitem{barvinokfixfirstcoeffs}
Alexander Barvinok, \emph{Computing the {E}hrhart quasi-polynomial of a
  rational simplex}, Math. Comp. \textbf{75} (2006), no.~255, 1449--1466
  (electronic), {\tt arXiv:math/0504444}.

\bibitem{nowhereharmonic}
Matthias Beck and Benjamin Braun, \emph{Nowhere-harmonic colorings of graphs},
  to appear in Proc. Amer. Math. Soc., {\tt arXiv:0907.1272}, 2011.

\bibitem{mahonianpartition}
Matthias Beck, Benjamin Braun, and Nguyen Le, \emph{Mahonian partition
  identities via polyhedral geometry}, to appear in Developments in
  Mathematics, {\tt arXiv:1103.1070}, 2011.

\bibitem{beckgesselleesavage}
Matthias Beck, Ira~M. Gessel, Sunyoung Lee, and Carla~D. Savage,
  \emph{Symmetrically constrained compositions}, Ramanujan J. \textbf{23}
  (2010), no.~1-3, 355--369, {\tt arXiv:0906.5573}.

\bibitem{ccd}
Matthias Beck and Sinai Robins, \emph{Computing the continuous discretely:
  Integer-point enumeration in polyhedra}, Undergraduate Texts in Mathematics,
  Springer, New York, 2007, Electronically available at {\tt
  http://math.sfsu.edu/beck/ccd.html}.

\bibitem{magiclatin}
Matthias Beck and Thomas Zaslavsky, \emph{An enumerative geometry for magic and
  magilatin labellings}, Ann. Comb. \textbf{10} (2006), no.~4, 395--413, {\tt
  arXiv:math.CO/0506315}.

\bibitem{iop}
\bysame, \emph{Inside-out polytopes}, Adv. Math. \textbf{205} (2006), no.~1,
  134--162, {\tt arXiv:math.CO/0309330}.

\bibitem{nnz}
\bysame, \emph{The number of nowhere-zero flows on graphs and signed graphs},
  J. Combin. Theory Ser. B \textbf{96} (2006), no.~6, 901--918, {\tt
  arXiv:math.CO/0309331}.

\bibitem{birkhoffcoloring}
George~D. Birkhoff, \emph{A determinant formula for the number of ways of
  coloring a map}, Ann. of Math. (2) \textbf{14} (1912/13), no.~1-4, 42--46.

\bibitem{breuerdall}
Felix Breuer and Aaron Dall, \emph{Bounds on the coefficients of tension and
  flow polynomials}, J. Algebraic Combin. \textbf{33} (2011), no.~3, 465--482,
  {\tt arXiv:1004.3470}.

\bibitem{breuersanyal}
Felix Breuer and Raman Sanyal, \emph{Ehrhart theory, modular flow reciprocity,
  and the {T}utte polynomial}, to appear in Math. Z., {\tt arXiv:0907.0845v1},
  2011.

\bibitem{deloerarambausantos}
Jes{\'u}s~A. De~Loera, J{\"o}rg Rambau, and Francisco Santos,
  \emph{Triangulations}, Algorithms and Computation in Mathematics, vol.~25,
  Springer-Verlag, Berlin, 2010.

\bibitem{ehrhartpolynomial}
Eug{\`e}ne Ehrhart, \emph{Sur les poly\`edres rationnels homoth\'etiques \`a
  {$n$}\ dimensions}, C. R. Acad. Sci. Paris \textbf{254} (1962), 616--618.

\bibitem{euler2}
Leonhard Euler, \emph{Demonstatio nonnullarum insignium proprietatum, quibus
  solida hedris planis inclusa sunt praedita}, Novi Comm. Acad. Sci. Imp.
  Petropol. \textbf{4} (1752/53), 140--160.

\bibitem{euler1}
\bysame, \emph{Elementa doctrinae solidorum}, Novi Comm. Acad. Sci. Imp.
  Petropol. \textbf{4} (1752/53), 109--140.

\bibitem{greeneacyclic}
Curtis Greene, \emph{Acyclic orientations}, Higher Combinatorics (M.~Aigner,
  ed.), NATO Adv. Study Inst. Ser., Ser. C: Math. Phys. Sci., vol.~31, Reidel,
  Dordrecht, 1977, pp.~65--68.

\bibitem{greenezaslavsky}
Curtis Greene and Thomas Zaslavsky, \emph{On the interpretation of {W}hitney
  numbers through arrangements of hyperplanes, zonotopes, non-{R}adon
  partitions, and orientations of graphs}, Trans. Amer. Math. Soc. \textbf{280}
  (1983), no.~1, 97--126.

\bibitem{haasenillpayne}
Christian Haase, Benjamin Nill, and Sam Payne, \emph{Cayley decompositions of
  lattice polytopes and upper bounds for {$h^*$}-polynomials}, J. Reine Angew.
  Math. \textbf{637} (2009), 207--216, {\tt arXiv:math/0804.3667}.

\bibitem{hibi}
Takayuki Hibi, \emph{Algebraic {C}ombinatorics on {C}onvex {P}olytopes},
  Carslaw, 1992.

\bibitem{linkerational}
Eva Linke, \emph{Rational {E}hrhart quasi-polynomials}, Preprint ({\tt
  arXiv:1006.5612v2}), 2011.

\bibitem{macdonald}
Ian~G. Macdonald, \emph{Polynomials associated with finite cell-complexes}, J.
  London Math. Soc. (2) \textbf{4} (1971), 181--192.

\bibitem{orlikterao}
Peter Orlik and Hiroaki Terao, \emph{Arrangements of hyperplanes}, Grundlehren
  der Mathematischen Wissenschaften [Fundamental Principles of Mathematical
  Sciences], vol. 300, Springer-Verlag, Berlin, 1992.

\bibitem{poincare}
Henri Poincar\'e, \emph{Sur la g\'en\'eralisation d'un theorem d'{E}uler
  relatif aux poly\`edres}, C. R. Acad. Sci. Paris (1893), 144--145.

\bibitem{samehrhartpolynomial}
Steven~V Sam, \emph{A bijective proof for a theorem of {E}hrhart}, Amer. Math.
  Monthly \textbf{116} (2009), no.~8, 688--701, {\tt arXiv:0801.4432v5}.

\bibitem{stanleythesis}
Richard~P. Stanley, \emph{Ordered structures and partitions}, American
  Mathematical Society, Providence, R.I., 1972, Memoirs of the American
  Mathematical Society, No. 119.

\bibitem{stanleyacyclic}
\bysame, \emph{Acyclic orientations of graphs}, Discrete Math. \textbf{5}
  (1973), 171--178.

\bibitem{stanleyec1}
\bysame, \emph{Enumerative {C}ombinatorics. {V}ol. 1}, Cambridge Studies in
  Advanced Mathematics, vol.~49, Cambridge University Press, Cambridge, 1997.

\bibitem{stanleyparkcityhyperplanes}
\bysame, \emph{An introduction to hyperplane arrangements}, Geometric
  combinatorics, IAS/Park City Math. Ser., vol.~13, Amer. Math. Soc.,
  Providence, RI, 2007, pp.~389--496.

\bibitem{stapledondelta}
Alan Stapledon, \emph{Inequalities and {E}hrhart {$\delta$}-vectors}, Trans.
  Amer. Math. Soc. \textbf{361} (2009), no.~10, 5615--5626, {\tt
  arXiv:math/0801.0873}.

\bibitem{stapledonadditive}
\bysame, \emph{Additive number theory and inequalities in {E}hrhart theory},
  Preprint ({\tt arXiv:0904.3035v2}), 2010.

\bibitem{whitneychromatic}
Hassler Whitney, \emph{A logical expansion in mathematics}, Bull. Amer. Math.
  Soc. \textbf{38} (1932), no.~8, 572--579.

\bibitem{zaslavskythesis}
Thomas Zaslavsky, \emph{Facing up to arrangements: face-count formulas for
  partitions of space by hyperplanes}, Mem. Amer. Math. Soc. \textbf{1} (1975),
  no.~154.

\bibitem{zaslavskybiasedgraphs7}
\bysame, \emph{Biased graphs. {VII}. {C}ontrabalance and antivoltages}, J.
  Combin. Theory Ser. B \textbf{97} (2007), no.~6, 1019--1040.

\bibitem{ziegler}
G{\"u}nter~M. Ziegler, \emph{Lectures on polytopes}, Springer-Verlag, New York,
  1995.

\end{thebibliography}

\setlength{\parskip}{0cm} 
\end{document}